\documentclass[portuges,12pt,letter]{article}
\usepackage[centertags]{amsmath}
\usepackage{amsfonts}
\usepackage{newlfont}
\usepackage{amscd}
\usepackage{graphics}
\usepackage{epsfig}
\usepackage{indentfirst}
\usepackage{amsxtra}
\usepackage[latin1]{inputenc}
\usepackage{amssymb, amsmath}
\usepackage{amsthm}
\usepackage[mathscr]{eucal}

\newtheorem{thm}{Theorem}[section]

\newtheorem{obe}[thm]{Remark}

\author{Fabio Silva Botelho \\ Department of Mathematics \\  Federal University of Santa Catarina, UFSC \\
Florian\'{o}polis, SC - Brazil}

\title{\bf  An approximate proximal numerical procedure concerning the generalized method of lines}
\date{}
\begin{document}
\maketitle

\abstract{This article develops an approximate proximal approach for the generalized method of lines. The present results are extensions and applications of previous ones which have been published since 2011, in books and articles such as \cite{901,909,120,700}. We also recall that in the generalized method of lines, the domain of the partial differential equation in question is discretized in lines (or in curves) and the concerning solution is developed on these lines, as functions of the boundary conditions and the domain boundary shape.}

\section{Introduction}
Let $\Omega \subset \mathbb{R}^2$ be an open, bounded and connected set where
$$\Omega=\{(x,y) \in \mathbb{R}^2\;:\;y_1(x) \leq y \leq y_2(x),\; a \leq x \leq b\}.$$

Here, we assume, $y_1,y_2:[a,b] \rightarrow \mathbb{R}$ are continuous functions.

Consider the Ginzburg-Landau type equation, defined by
\begin{equation}\left\{\begin{array}{ll}
-\varepsilon \nabla^2 u+\alpha u^3-\beta u=f,& \text{ in } \Omega,
 \\
u=0, & \text{ on } \partial \Omega. \end{array} \right.\end{equation}

Here $\varepsilon>0$, $\alpha>0,\; \beta>0$ and $f \in L^2(\Omega)$.

Also, $u \in W_0^{1,2}(\Omega)$ and the equation in question must be considered in a distributional sense.

\section{ The numerical method}
We discretize the interval $[a,b]$ into $N$ same measure sub-intervals, through a partition $$P=\{x_0=a,x_1,\cdots, x_N=b\},$$
where $x_n=a+nd,\; \forall n \in \{1,\cdots,N-1\}.$ Here $$d=\frac{(b-a)}{N}.$$

Through such a procedure, we generate $N$ vertical lines parallel to the Cartesian axis  $0y$, so that for each line $n$ based on the point $x_n$
we are going to compute an approximate solution $u_n(y)$ corresponding to values of $u$ on such a line.

Considering this procedure, the equation system obtained in partial finite differences is given by

$$-\varepsilon\left( \frac{u_{n+1}-2u_n+u_{n-1}}{d^2}+ \frac{\partial^2 u_n}{\partial y^2}\right)+\alpha u_n^3-\beta u_n=f_n,$$
$\forall n \in \{1, \cdots,N-1\},$ with the boundary conditions $$u_0=0,$$ and $$u_N=0.$$

Let $K>0$ be an appropriate constant to be specified.

In a proximal approach, considering an initial solution $$\{(u_0)_n\}$$ we redefine the system of equations in question as below indicated.
$$-\varepsilon \left(\frac{u_{n+1}-2u_n+u_{n-1}}{d^2}+ \frac{\partial^2 u_n}{\partial y^2}\right)+\alpha u_n^3-\beta u_n +K u_n-K(u_0)_n=f_n,$$
$\forall n \in \{1, \cdots,N-1\},$ with the boundary conditions $$u_0=0,$$ and $$u_N=0.$$

Hence, we may denote
$$u_{n+1}-\left(2+K\frac{d^2}{\varepsilon}\right)u_n+u_{n-1}+T(u_n)+ \tilde{f}_n\frac{d^2}{\varepsilon}=0,$$
where
$$ T(u_n)=\left(-\alpha u_n^3+\beta u_n\right)\frac{d^2}{\varepsilon}+\frac{\partial^2 u_n}{\partial y^2}d^2,$$
and $\tilde{f}_n=K(u_0)_n+f_n,$ $\forall n \in \{1,\cdots,N-1\}.$

In particular, for $n=1$, we get

$$u_{2}-\left(2+K\frac{d^2}{\varepsilon}\right)u_1+T(u_1)+ \tilde{f}_1\frac{d^2}{\varepsilon}=0,$$
so that
$$u_1=a_1u_2+b_1T(u_2)+c_1+E_1,$$
where
$$a_1=\frac{1}{ 2+K \frac{d^2}{\varepsilon}},$$
$$b_1=a_1,$$
and $$c_1=a_1 \tilde{f_1}\frac{d^2}{\varepsilon}$$
and the error $E_1$, proportional to $1/K$, is given by
$$E_1=b_1(T(u_1)-T(u_2)).$$

Now, reasoning inductively, having
$$u_{n-1}=a_{n-1}u_n+b_{n-1}T(u_n)+c_{n-1}+E_{n-1}$$
for the line $n$, we have
\begin{eqnarray}
&&u_{n+1}-\left(2+K\frac{d^2}{\varepsilon}\right)u_n+a_{n-1}u_n+b_{n-1}T(u_n)+c_{n-1}+E_{n-1}
\nonumber \\ &&+T(u_n)+ \tilde{f}_n\frac{d^2}{\varepsilon}=0,\end{eqnarray}
so that
$$u_n=a_nu_{n+1}+b_nT(u_n)+c_n+E_n,$$
where
$$a_n=\frac{1}{ 2+K \frac{d^2}{\varepsilon}-a_{n-1}},$$
$$b_{n}=a_n(b_{n-1}+1),$$
and $$c_n=a_n\left(c_{n-1}+ \tilde{f}_n\frac{d^2}{\varepsilon}\right)$$
and the error $E_n$, is given by
$$E_n=a_nE_{n-1}+b_n(T(u_n)-T(u_{n+1})),$$
$\forall n \in \{1,\cdots,N-1\}.$

In particular, for $n=N-1$, we have $u_N=0$ so that,
\begin{eqnarray}
u_{N-1} &\approx& a_{N-1}u_N+b_{N-1} T(u_N)+c_{N-1} \nonumber \\ &\approx&
a_{N-1}u_N+b_{N-1}\frac{\partial^2 u_{N-1}}{\partial y^2} d^2+b_{N-1}(-\alpha u_N^3+\beta u_N)\frac{d^2}{\varepsilon}+c_{N-1} \nonumber \\ &=&
b_{N-1}\frac{\partial^2 u_{N-1}}{\partial y^2} d^2+c_{N-1}.
\end{eqnarray}

This last equation is an ODE from which we may easily obtain $u_{N-1}$ with the boundary conditions $$u_{N-1}(y_1(x_{N-1}))=u_{N-2}(y_2(x_{N-1}))=0.$$

Having $u_{N-1}$, we may obtain $u_{N-2}$  though the equation
\begin{eqnarray}
u_{N-2} &\approx& a_{N-2}u_{N-1}+b_{N-2} T(u_{N-1})+c_{N-2} \nonumber \\ &\approx&
a_{N-2}u_{N-1}+b_{N-2}\frac{\partial^2 u_{N-2}}{\partial y^2} d^2+b_{N-2}(-\alpha u_{N-1}^3+\beta u_{N-1})\frac{d^2}{\varepsilon}+c_{N-2},
\end{eqnarray}
with the boundary conditions $$u_{N-2}(y_1(x_{N-2}))=u_{N-2}(y_2(x_{N-2}))=0.$$
An so on, up to finding $u_1.$

The next step is to replace $\{(u_0)_n\}$ by $\{u_n\}$ and then to repeat the process until an appropriate convergence criterion is satisfied.

The problem is then approximately solved.
\section{ A numerical example}
We present numerical results for $\Omega=[0,1] \times [0,1],$ $ \alpha=\beta=1$, $f\equiv 1$  in $\Omega$, $N=100$, $K=50$ and for
$$\varepsilon=0.1, \; 0.01 \text{ and } 0.001.$$

\begin{figure}
\centering \includegraphics [width=4in]{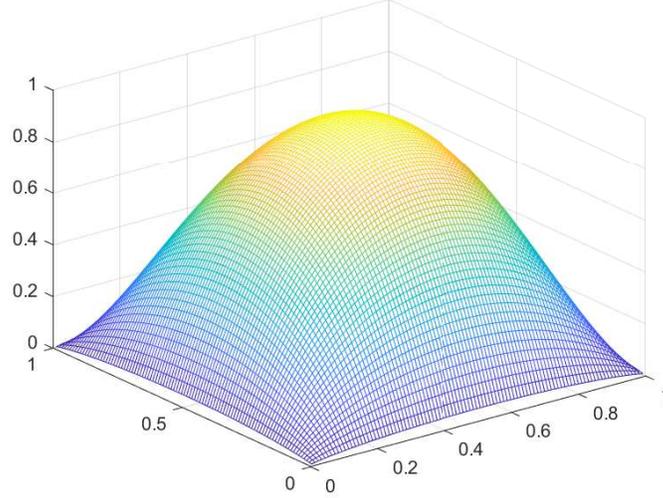}
\caption{\small{Solution u for $\varepsilon=0.1$}\label{figLGMay2022-3}}
\end{figure}

\begin{figure}
\centering \includegraphics [width=4in]{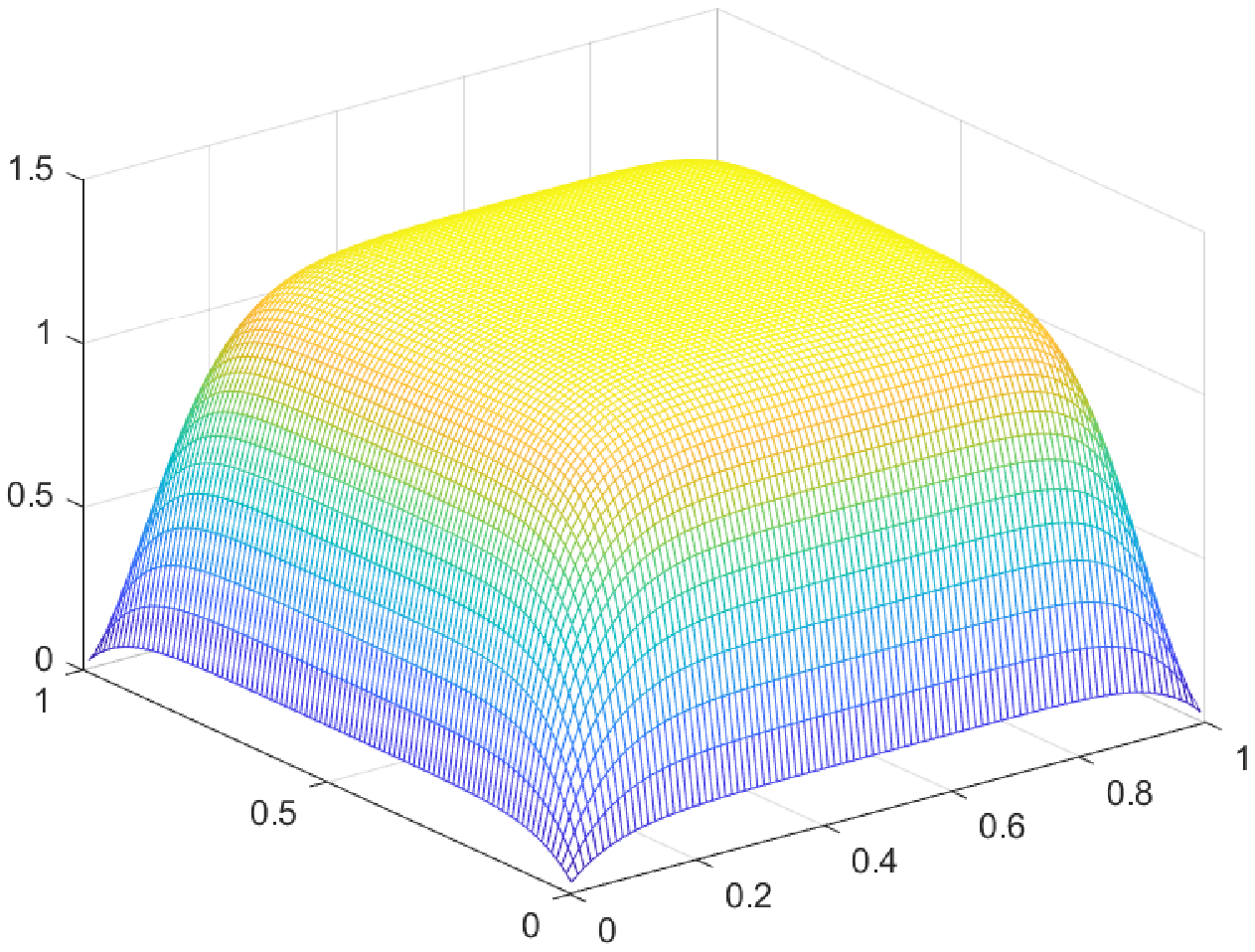}
\caption{\small{Solution u for $\varepsilon=0.01$}\label{figLGMay2022-2}}
\end{figure}

\begin{figure}
\centering \includegraphics [width=4in]{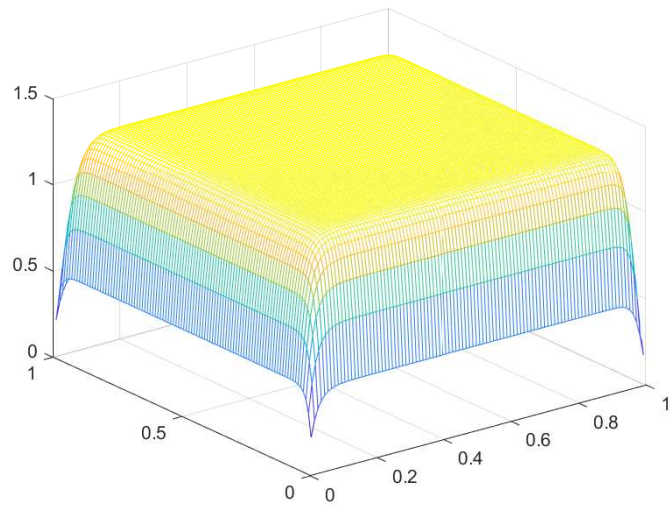}
\caption{\small{Solution u for $\varepsilon=0.001$}\label{figLGMay2022-1}}
\end{figure}

For such values of $\varepsilon$, please see figures \ref{figLGMay2022-3}, \ref{figLGMay2022-2} and  \ref{figLGMay2022-1}, respectively.

\begin{obe} We observe that as $\varepsilon>0$ decreases to the value $0.001$ the solution approaches the constant value $1.3247$ along the domain, up to the satisfaction of  boundary conditions. This is expect, since this value is an approximate solution of equation $u^3-u-1=0.$
\end{obe}

\section{A general proximal explicit approach}

Based on the algorithm presented in the last section, we develop a software in MATHEMATICA in order to approximately solve the following
equation.
\begin{equation}\left\{\begin{array}{ll}
-\varepsilon\left(\frac{\partial^2 u}{\partial r^2}+\frac{1}{r}\frac{\partial u}{\partial r}+\frac{1}{r^2} \frac{\partial^2 u}{\partial \theta^2}\right)  u+\alpha u^3-\beta u=f,& \text{ in } \Omega,
 \\
u=0, & \text{ on } \partial \Omega_1, \\
u=u_f(\theta), & \text{ on } \partial \Gamma_2. \end{array} \right.\end{equation}

Here $$\Omega=\{(r,\theta) \in \mathbb{R}^2\;:\; 1\leq r \leq 2,\; 0\leq \theta \leq 2\pi\},$$
$$\partial \Omega_1=\{(1,\theta) \in \mathbb{R}^2\;:\; 0\leq \theta \leq 2\pi\},$$
$$\partial \Omega_2=\{(2,\theta) \in \mathbb{R}^2\;:\; 0\leq \theta \leq 2\pi\},$$
$\alpha=\beta=1$, $K=10$, $N=100$ and $f\equiv 1, \text{ on } \Omega.$
\\
\\
At this point we present such a software in MATHEMATICA. 
\\
\\
******************************************
\begin{enumerate}
\item $m8 = 100;$
\item $d = 1.0/m8;$
\item $K = 10.0;$
\item $e1 = 0.01;$
\item $A = 1.0;$
\item $B = 1.0;$

\item For[i = 1, i $<$ m8 + 1, i++,

  uo[i] = 0.0];
  
\item For[k = 1, k $<$ 150, k++,

 Print[k];
 
 a[1] = 1/(2.0 + $K*d^2/e1$);
 
 b[1] = a[1];
 
 c[1] = a[1]*(K*uo[1] + 1.0)*$d^2/e1$;
 
 For[i = 2, i $<$ m8, i++,
 
  a[i] = $1/(2.0 + K*d^2/e1 - a[i - 1]);$
  
  b[i] = a[i]*(b[i - 1] + 1);
  
  c[i] = a[i]*(c[i - 1] + $(K*uo[i] + 1.0)*d^2/e1)$];
  
 \item u[m8] = uf[x]; d1 = 1.0;
 
 \item For[i = 1, i $<$ m8, i++,
 
  t[m8 - i] = 1 + (m8 - i)*d;
  
  A1 = (a[m8 - i]*u[m8 - i + 1] +
  
      b[m8 - i]*(-A*$u[m8 - i + 1]^3$ + B*u[m8 - i + 1])*$d^2/e1*d1^2$ +
      
      c[m8 - i] +
      
      $d^2*d1^2$*b[m8 - i]*(D[u[m8 - i + 1], {x, 2}]/$t[m8 - i]^2$) +
      
     $ d1^2*1/t[m8 - i]$*b[m8 - i]*
       $d^2$ (uo[m8 - i + 1] - uo[m8 - i])/d)/(1.0);
      
  A1 = Expand[A1];
  
  A1 = Series[
    A1, \{uf[x], 0, 3\}, \{uf'[x], 0, 1\}, \{uf''[x], 0, 1\}, \{uf'''[x], 0,
     0\}, \{uf''''[x], 0, 0\}];
     
  A1 = Normal[A1];
  
  u[m8 - i] = Expand[A1]];
  
 For[i = 1, i $<$ m8 + 1, i++,
 
  uo[i] = u[i]]; d1 = 1.0;
  
 Print[Expand[u[m8/2]]]]
 \end{enumerate}
 
 ************************************
\\
\\
For such a general approach, for $\varepsilon=0.1$, we have obtained the following lines (here $x$ stands for $\theta.$).

\begin{eqnarray}
u[10](x)&=&0.4780 +0.0122\; u_f[x]-0.0115\; u_f[x]^2+0.0083\; u_f[x]^3+0.00069 u_f''[x]
\nonumber \\ &&-0.0014\; u_f[x] (u_f'')[x]+0.0016\; uf[x]^2 u_f'[x]-0.00092\; u_f[x]^3 u_f''[x]
\end{eqnarray}

\begin{eqnarray}u[20](x)&=&0.7919 +0.0241\; u_f[x]-0.0225\; u_f[x]^2+0.0163 \;u_f[x]^3+0.0012\; u_f''[x]
\nonumber \\ &&-0.0025\; u_f[x] (u_f'')[x]+0.0030\; u_f[x]^2 (u_f'')[x]-0.0018\; u_f[x]^3 (u_f'')[x].\end{eqnarray}

\begin{eqnarray}u[30](x)&=&0.9823 +0.0404\; u_f[x]-0.0375\; u_f[x]^2+0.0266\; u_f[x]^3+0.00180\; (u_f'')[x]
\nonumber \\ &&-0.00362 \;u_f[x] (u_f'')[x]+0.0043\; u_f[x]^2 (uf'')[x]-0.0028 \;u_f[x]^3 (uf'')[x]
\end{eqnarray}

\begin{eqnarray}u[40](x)&=&1.0888 +0.0698\; u_f[x]-0.0632\; u_f[x]^2+0.0433\; u_f[x]^3+0.0026 (u_f'')[x]
\nonumber \\ &&-0.0051\; u_f[x] (u_f'')[x]+0.0061 u_f[x]^2 (u_f'')[x]-0.0043 u_f[x]^3 (u_f'')[x] \end{eqnarray}

\begin{eqnarray}u[50](x)&=&1.1316 +0.1277\; u_f[x]-0.1101 \;u_f[x]^2+0.0695 \;u_f[x]^3+0.0037\; (u_f'')[x]
\nonumber \\ &&-0.0073 \;u_f[x] (u_f'')[x]+0.0084\; u_f[x]^2 (u_f'')[x]-0.0062\; u_f[x]^3 (u_f'')[x] \end{eqnarray}

\begin{eqnarray}u[60](x)&=&1.1104 +0.2389\; u_f[x]-0.1866\; u_f[x]^2+0.0988\; u_f[x]^3+0.0053 (u_f'')[x]
\nonumber \\ &&-0.0099\; u_f[x] (u_f'')[x]+0.0105\; u_f[x]^2 (u_f'')[x]-0.0075\; u_f[x]^3 (u_f''])[x]
\end{eqnarray}

\begin{eqnarray}u[70](x)&=&1.0050 +0.4298\; u_f[x]-0.273813\; u_f[x]^2+0.0949\; u_f[x]^3+0.0070\; (u_f'')[x]
\nonumber \\ &&-0.0116\; u_f[x] (u_f'')[x]+0.0102\; u_f[x]^2 (u_f'')[x]-0.0061\; u_f[x]^3 (u_f'')[x]
\end{eqnarray}

\begin{eqnarray}u[80](x)&=&0.7838 +0.6855\; u_f[x]-0.2892\; u_f[x]^2+0.0161\; u_f[x]^3+0.0075\; (u_f'')[x]
\nonumber \\ &&-0.0098 u_f[x] (u_f'')[x]+0.0063084 uf[x]^2 (uf'')[x]-0.0027 u_f[x]^3 (u_f'')[x]
\end{eqnarray}

\begin{eqnarray}u[90](x)&=&0.4359 +0.9077\; u_f[x]-0.1621 u_f[x]^2-0.0563\; u_f[x]^3+0.0051\; (u_f'')[x]
\nonumber \\ &&-0.0047 \;u_f[x] (u_f'')[x]+0.0023\; u_f[x]^2 (u_f'')[x]-0.00098\; u_f[x]^3 (u_f'')[x]
\end{eqnarray}

For $\varepsilon=0.01$, we have obtained the following line expressions.

\begin{eqnarray}
u[10](x)&=&1.0057 +2.07*10^{-11}\; u_f[x]-1.85*10^{-11}\; u_f[x]^2+1.13*10^{-11}\;u_f[x]^3\nonumber \\ &&+4.70*10^{-13} \;(u_f'')[x]
-8.44*10^{-13} u_f[x] (u_f'')[x]\nonumber \\ &&+7.85*10^{-13} \;u_f[x]^2 (u_f'')[x]-6.96*10^{-14}\; u_f[x]^3 (u_f'')[x]
\end{eqnarray}

\begin{eqnarray}
u[20](x)&=&1.2512 +2.13*10^{-10}\; u_f[x]-1.90*10^{-10}\; u_f[x]^2+1.16*10^{-10}\; u_f[x]^3\nonumber \\ &&+3.94*10^{-12}\; (u_f'')[x]
-7.09*10^{-12} u_f[x] (uf'')[x]+6.61*10^{-12}\nonumber \\ &&\; u_f[x]^2 (u_f'')[x]-7.17*10^{-13}\; u_f[x]^3 (u_f'')[x]
\end{eqnarray}

\begin{eqnarray}u[30](x)&=&1.3078 +3.80*10^{-9} uf[x]-3.39*10^{-9}\; u_f[x]^2+2.07*10^{-9}\; u_f[x]^3\nonumber \\ &&+5.65*10^{-11} (u_f'')[x]
-1.018*10^{-10}\; u_f[x] (u_f'')[x]\nonumber \\ &&+9.52*10^{-11}\; u_f[x]^2 (u_f'')[x]-1.27*10^{-11}\; uf[x]^3 (u_f'')[x]
\end{eqnarray}

\begin{eqnarray}
u[40](x)&=&1.3208 +7.82*10^{-8} u_f[x]-6.98*10^{-8}\; u_f[x]^2+4.27*10^{-8}\; u_f[x]^3\nonumber \\ &&+9.27*10^{-10}\; (u_f'')[x]
-1.67*10^{-9} u_f[x] (u_f'')[x]\nonumber \\ &&+1.57*10^-9 \;u_f[x]^2 (u_f'')[x]-2.62*10^{-10}\; u_f[x]^3 (u_f'')[x]
\end{eqnarray}

\begin{eqnarray}
u[50](x)&=&1.3238 +1.67*10^{-6}\; u_f[x]-1.49*10^{-6}\; u_f[x]^2+9.15*10^{-7}\; uf[x]^3\nonumber \\ &&+1.54*10^{-8}\; (u_f'')[x]
-2.79*10^{-8}\; u_f[x] (u_f'')[x]\nonumber \\ &&+2.64*10^{-8}\; u_f[x]^2 (u_f'')[x]-5.62*10^{-9}\; u_f[x]^3 (u_f'')[x]
\end{eqnarray}

\begin{eqnarray}
u[60](x)&=& 1.32449 +0.000036\; u_f[x]-0.000032\; u_f[x]^2+0.000019\; u_f[x]^3\nonumber \\ &&+2.51*10^{-7} \;(u_f'')[x]
-4.57*10^{-7}\; u_f[x] (u_f'')[x]\nonumber \\ &&+4.36*10^{-7}\; u_f[x]^2 (u_f'')[x]-1.21*10^{-7}\; u_f[x]^3 (u_f'')[x]
\end{eqnarray}

\begin{eqnarray}
u[70](x)&=&1.32425 +0.00079\; u_f[x]-0.00070\; u_f[x]^2+0.00043\; u_f[x]^3\nonumber \\ &&+3.89*10^{-6}\; (u_f'')[x]
-7.12*10^{-6}\; u_f[x] (u_f'')[x]\nonumber \\ &&+6.89*10^{-6}\; u_f[x]^2 (u_f'')[x]-2.64*10^{-6} \;u_f[x]^3 (u_f'')[x]
\end{eqnarray}

\begin{eqnarray}u[80](x)&=&1.31561 +0.017\; u_f[x]-0.015\; u_f[x]^2+0.009\; u_f[x]^3\nonumber \\ &&+0.000053\; (u_f'')[x]
-0.000098\; u_f[x] (u_f '')[x]\nonumber \\ &&+0.000095\; u_f[x]^2 (u_f'')[x]-0.000051\; u_f[x]^3 (u_f'')[x]
\end{eqnarray}

\begin{eqnarray}u[90](x)&=&1.14766 +0.296 u_f[x]-0.1991\; u_f[x]^2+0.0638\; u_f[x]^3\nonumber \\ &&+0.00044\; (u_f'')[x]
-0.00067\; u_f[x] (u_f'')[x]\nonumber \\ &&+0.00046\; u_f[x]^2 (u_f'')[x]-0.00018\; u_f[x]^3 (u_f'')[x]
\end{eqnarray}

\begin{obe} We observe that as $\varepsilon>0$ decreases to the value $0.01$ the solution approaches the constant value $1.3247$ along the domain, up to the satisfaction of  boundary conditions. This is expect, since this value is an approximate solution of equation $u^3-u-1=0.$
\end{obe}

\end{document}